\documentclass[12pt]{amsart}

\usepackage[paperheight=11in,paperwidth=8.5in,left=1in,top=1.25in,right=1in,bottom=1.25in]{geometry}
\usepackage[pdfstartpage={1},pdfstartview={FitH},pdftitle={},pdfauthor={},pdfsubject={},pdfcreator={},pdfproducer={},pdfkeywords={},pagebackref=false,colorlinks=true,linkcolor=black,citecolor=black,filecolor=black,urlcolor=black,]{hyperref}
\usepackage{afterpage}
\usepackage{amsfonts}
\usepackage{amsmath}
\allowdisplaybreaks[4]
\usepackage{amssymb}
\usepackage{amsthm}
\usepackage{array}
\usepackage{bbm}
\usepackage[justification=centering]{caption}
\usepackage[capitalize,nameinlink,noabbrev,nosort]{cleveref}
\usepackage{enumerate}
\usepackage{float}
\restylefloat{figure}
\usepackage{graphicx}
\usepackage{mathrsfs} 
\usepackage{mathtools}
\usepackage[framemethod=tikz,ntheorem,xcolor]{mdframed}
\usepackage{parcolumns}
\usepackage{tabularx}
\usepackage{tikz}\pdfpageattr{/Group <</S /Transparency /I true /CS /DeviceRGB>>}
\usetikzlibrary{arrows,calc,intersections,matrix,positioning,through,patterns}
\usepackage{tikz-cd}
\tikzset{commutative diagrams/.cd,every label/.append style = {font = \normalsize}}
\usepackage[all]{xy}
\usepackage[aligntableaux=center]{ytableau}

\usetikzlibrary{external}
\newtheorem{thm}{Theorem}
\newtheorem{lem}{Lemma}
\theoremstyle{definition}
\newtheorem{rmk}{Remark}
\theoremstyle{definition}
\newtheorem{example}{Example}
\theoremstyle{definition}
\newtheorem{defn}{Definition}
\theoremstyle{remark}

\AfterEndEnvironment{claim}{\noindent\ignorespaces}
\newtheorem*{claimpf_no_qed}{Proof of Claim}

\AfterEndEnvironment{claimpf}{\noindent\ignorespaces}

\newcommand{\rf}[1]{\hyperref[#1]{(\ref*{#1})}}

\DeclareMathOperator{\SL}{SL}

\DeclareMathOperator{\supp}{supp}

\def\RR{{\mathbb R}}

\newcommand{\exps}[1]{a(#1)}
\def\expts{\exps t}

\def\Z{{\mathbb Z}}

\def\gfr{{\mathfrak g}}

\def\P{\mathcal{P}}
\def\Ptnn{\P_{\geq0}}
\def\Ptp{\P_{>0}}

\def\dotP{{\dot\P}}
\def\dotPtnn{\dotP_{\geq0}}
\def\dotPtp{\dotP_{>0}}

\def\B{\mathcal{B}} 
\def\Btnn{\B_{\geq0}}
\def\Btp{\B_{>0}}

\def\dotB{{\dot\B}}
\def\dotBtnn{\dotB_{\geq0}}
\def\dotBtp{\dotB_{>0}}

\def\Gtnn{G_{\geq0}}
\def\Gtp{G_{>0}}

\def\dotG{{\dot G}}
\def\dotGtnn{\dotG_{\geq0}}
\def\dotGtp{\dotG_{>0}}

\def\Utp{U_{>0}}
\def\Ttp{T_{>0}}

\def\ProjSpace{\mathbb{P}}
\def\PSl{\ProjSpace(\irrep_\l)}
\def\PSltnn{\PSl_{\geq0}}
\def\PSltp{\PSl_{>0}}

\def\l{\lambda}

\def\<{\langle}
\def\>{\rangle}

\newcommand{\CB}[1]{{ _{#1}\mathbf{B}}}
\def\CBl{\CB{\l}}

\def\DomWeights{X^+}
\def\irrep{\Lambda}

\def\line{L}

\def\dotT{{\dot T}}
\def\dotBB{{\dot B}}
\def\dotx{{\dot x}}
\def\doty{{\dot y}}
\def\dotI{{\dot I}}
\def\dotJ{{\dot J}}
\def\barI{{I}}
\def\barJ{{J}}

\newcommand{\LusRef}[1]{\S#1}
\newcommand{\Lus}[1]{\cite[\LusRef{#1}]{Lus}}
\newcommand{\LusGP}[1]{\cite[\LusRef{#1}]{LusGP}}

\title{The totally nonnegative part of $G/P$ is a ball}

\author{Pavel Galashin}
\address{Department of Mathematics, Massachusetts Institute of Technology, 77 Massachusetts Avenue,
Cambridge, MA 02139, USA}
\email{\href{mailto:galashin@mit.edu}{galashin@mit.edu}}
\author{Steven N. Karp}
\address{Department of Mathematics, University of Michigan, 2074 East Hall, 530 Church Street, Ann Arbor, MI 48109-1043, USA}
\email{\href{mailto:snkarp@umich.edu}{snkarp@umich.edu}}
\email{\href{mailto:tfylam@umich.edu}{tfylam@umich.edu}}
\author{Thomas Lam}
\thanks{T.L.\ acknowledges support from the NSF under agreement No.\ DMS-1464693.}

\begin{document}

\begin{abstract}
We show that the totally nonnegative part of a partial flag variety (in the sense of Lusztig) is homeomorphic to a closed ball.
\end{abstract}

\date{\today}

\maketitle

\numberwithin{equation}{section}
\section{Introduction}
\noindent Let $G$ be a simply connected, semisimple algebraic group defined and split over $\RR$, and let $\P^\barJ$ be the (real) partial flag variety of $G$ associated to a subset $\barJ$ of the set $\barI$ of simple roots.  Lusztig \cite{Lus} has defined a ``remarkable polyhedral subspace'' $\P^\barJ_{\geq 0}$ of $\P^\barJ$, called the \emph{totally nonnegative part}.  In this paper we establish the following result.

\begin{thm}\label{main}
The totally nonnegative part $\P^\barJ_{\ge 0}$ is homeomorphic to a closed ball of dimension equal to that of $\P^\barJ$.
\end{thm}

Lusztig proved that $\P^\barJ_{\ge 0}$ is contractible~\cite[\S 4.4]{LusIntro}. He also defined a partition of $\P^\barJ_{\ge 0}$~\Lus{8.15}, which was shown to be a cell decomposition by Rietsch~\cite{Rietsch}. Williams conjectured that this cell decomposition forms a regular CW complex~\cite[\S 7]{Wil}. Rietsch and Williams showed that $\P^\barJ_{\ge 0}$ forms a CW complex~\cite{RW-CW}, and have proved Williams's conjecture up to homotopy~\cite{RW}. These results complement work by Hersh \cite{Hersh} in the unipotent group case.

Our proof of \cref{main} employs the vector field $\tau$ on $\P^\barJ$ which is the sum of all Chevalley generators of $\gfr$, which Lusztig used to show that $\P^\barJ_{\ge 0}$ is contractible. The flow defined by $\tau$ is a {\itshape contractive flow} on $\P^\barJ_{\ge 0}$ in the sense of our paper~\cite{GKL}, so in particular it contracts all of $\P^\barJ_{\ge 0}$ to a unique fixed point $p_0\in\P^\barJ_{\ge 0}$. The machinery of~\cite{GKL} (see \cref{lem:top}) generates a homeomorphism from $\P^\barJ_{\ge 0}$ to a closed ball $B\subset \P^\barJ_{\ge0}$ centered at $p_0$, by mapping each trajectory in $\P^\barJ_{\ge 0}$ to its intersection with $B$. This generalizes~\cite[Theorem 1.1]{GKL} from type $A$ Grassmannians to all partial flag varieties. We remark that \cref{main} is new in all other cases, including for the complete flag variety and multi-step flag varieties in $\RR^n$.

\section{Preliminaries}

\noindent In this section, we recall some background from Lusztig~\cite{Lus,LusGP} and~\cite{GKL}.

\subsection{Pinnings}

Let $\gfr$ denote the Lie algebra of $G$ over $\RR$. We fix Chevalley generators $(e_i,f_i)_{i\in \barI}$ of $\gfr$, so that the elements $h_i:=[e_i,f_i]$ ($i\in\barI$) span the Lie algebra of a split real maximal torus $T$ of $G$. For $i\in\barI$ and $t\in\mathbb{R}$, we define the elements of $G$
\[
x_i(t) := \exp(te_i),\quad y_i(t) := \exp(tf_i).
\]
We also let $\alpha_i^\vee:\mathbb{R}^\ast\to T$ be the homomorphism of algebraic groups whose tangent map takes $1\in\mathbb{R}$ to $h_i$. The $x_i(t)$'s (respectively, $y_i(t)$'s) generate the unipotent radical $U^+$ of a Borel subgroup $B^+$ (respectively, $U^-$ and $B^-$) of $G$, with $B^+\cap B^- = T$. The data $(T,B^+,B^-,x_i,y_i;i\in \barI)$ is called a \emph{pinning} for $G$.

\subsection{Total positivity}

Let $i_1,\dots,i_\ell$ be a sequence of elements of $\barI$ such that the product of simple reflections $s_{i_1}\cdots s_{i_\ell}$ is a reduced decomposition of the longest element of the Weyl group. We define the totally positive parts
\begin{align*}
\Utp^+:=\{x_{i_1}(t_1)\cdots x_{i_\ell}(t_\ell) : t_1,\dots,t_\ell>0 \},\quad
\Utp^-:=\{y_{i_1}(t_1)\cdots y_{i_\ell}(t_\ell) : t_1,\dots,t_\ell>0 \}.
\end{align*}
Let $\Gtp:=\Utp^+ \Ttp\Utp^-=\Utp^- \Ttp\Utp^+$, where $\Ttp$ is generated by $\alpha_i^\vee(t)$ for $i\in \barI$ and $t>0$.

The \emph{complete flag variety} $\B$ of $G$ is the space of all Borel subgroups $B$ of $G$. Its totally positive part is
\[
\Btp:=\{u B^+ u^{-1} : u\in \Utp^-\}=\{u B^- u^{-1} : u\in \Utp^+\},
\]
and its totally nonnegative part $\Btnn$ is the closure of $\Btp$ in $\B$.

We now fix a subset $\barJ\subset \barI$ and define $P_\barJ^+$ to be the subgroup of $G$ generated by $B^+$ and $\{y_j(t) : j\in \barJ,\, t\in \RR\}$. Let $\P^\barJ:=\{g P_\barJ^+ g^{-1} : g\in G\}$ be a \emph{partial flag variety} of $G$. (In the case $J=\emptyset$, we have $\P^\barJ = \B$.) For each $B\in\B$, there is a unique $P\in\P^\barJ$ such that $B\subset P$; denote by $\pi^\barJ:\B\to\P^\barJ$ the natural projection map that sends $B\in\B$ to the $P\in\P^\barJ$ that contains it. We define the totally positive and totally nonnegative parts
\begin{align}\label{eq:pi_J}
\Ptp^\barJ:=\pi^\barJ(\Btp),\quad\Ptnn^\barJ:=\pi^\barJ(\Btnn).
\end{align}

\begin{example}\label{ex:flag}
Let $G:=\SL_3(\RR)$. We may take a pinning with $I = \{1, 2\}$ and
\begin{gather*}
e_1 = \scalebox{0.80}{$\begin{bmatrix}0 & 1 & 0 \\ 0 & 0 & 0 \\ 0 & 0 & 0\end{bmatrix}$},\;
f_1 = \scalebox{0.80}{$\begin{bmatrix}0 & 0 & 0 \\ 1 & 0 & 0 \\ 0 & 0 & 0\end{bmatrix}$},\;
x_1(t) = \scalebox{0.80}{$\begin{bmatrix}1 & t & 0 \\ 0 & 1 & 0 \\ 0 & 0 & 1\end{bmatrix}$},\;
y_1(t) = \scalebox{0.80}{$\begin{bmatrix}1 & 0 & 0 \\ t & 1 & 0 \\ 0 & 0 & 1\end{bmatrix}$},\;
\alpha_1^\vee(t) = \scalebox{0.80}{$\begin{bmatrix}t & 0 & 0 \\ 0 & t^{-1} & 0 \\ 0 & 0 & 1\end{bmatrix}$},\\[2pt]
e_2 = \scalebox{0.80}{$\begin{bmatrix}0 & 0 & 0 \\ 0 & 0 & 1 \\ 0 & 0 & 0\end{bmatrix}$},\;
f_2 = \scalebox{0.80}{$\begin{bmatrix}0 & 0 & 0 \\ 0 & 0 & 0 \\ 0 & 1 & 0\end{bmatrix}$},\;
x_2(t) = \scalebox{0.80}{$\begin{bmatrix}1 & 0 & 0 \\ 0 & 1 & t \\ 0 & 0 & 1\end{bmatrix}$},\;
y_2(t) = \scalebox{0.80}{$\begin{bmatrix}1 & 0 & 0 \\ 0 & 1 & 0 \\ 0 & t & 1\end{bmatrix}$},\;
\alpha_2^\vee(t) = \scalebox{0.80}{$\begin{bmatrix}1 & 0 & 0 \\ 0 & t & 0 \\ 0 & 0 & t^{-1}\end{bmatrix}$}.
\end{gather*}
Then $B^+$, $B^-$, and $T$ are the subgroups of upper-triangular, lower-triangular, and diagonal matrices, respectively. Taking $J:=\emptyset$, we can identify $\P^\barJ=\B$ with the space of complete flags in $\RR^3$. Explicitly, if $g = [g_1 \,|\, g_2\,|\, g_3]\in G$, then $gB^+g^{-1}\in\P^\barJ$ corresponds to the flag $\{0\} \subset \langle g_1\rangle \subset \langle g_1, g_2\rangle \subset \RR^3$. One can check that the totally nonnegative part $\Ptnn^\barJ$ consists of flags $\{0\} \subset V_1 \subset V_2 \subset \RR^3$ where $V_1$ is spanned by a vector $(v_1, v_2, v_3)\in\RR^3$ with $v_1, v_2, v_3\ge 0$, and $V_2$ is orthogonal to a vector $(w_1, -w_2, w_3)\in\RR^3$ with $w_1, w_2, w_3 \ge 0$. Hence we may identify $\Ptnn^\barJ$ with the subset of $\RR^6$ of points $(v_1, v_2, v_3, w_1, w_2, w_3)$ satisfying the (in)equalities
\[
v_1 + v_2 + v_3 = 1,\;\;\; w_1 + w_2 + w_3 = 1,\;\;\; v_1w_1 - v_2w_2 + v_3w_3 = 0,\;\;\; v_1, v_2, v_3, w_1, w_2, w_3 \ge 0.
\]
\cref{main} implies that this region is homeomorphic to a $3$-dimensional closed ball. Its cell decomposition (determined by which of the six coordinates are zero) is shown in \cref{fig:pokeball}.
\end{example}

  \newcommand\pbase{(3cm and 1cm)}
  \newcommand\pone{(1cm and 3cm)}
  \newcommand\ptwo{(0.6cm and 3cm)}
  \newcommand\pcirc{(3cm and 3cm)}

  \def\lw{2pt}
  \def\opF{0.6}
  \def\opFF{1}
  
  \def\opB{0.5}
  \def\opBB{0.2}

  \begin{figure}

\begin{tikzpicture}[scale=0.86]
\useasboundingbox(-4.6cm,-3.2cm)rectangle(4.6cm,3.2cm);
        \path[name path=base] (0,0) ellipse \pbase;
        \path[name path=one] (0,0) ellipse \pone;
        \path[name path=two] (0,0) ellipse \ptwo;
        \path[name path=circ] (0,0) ellipse \pcirc;
        
        \path [name intersections={of=base and one}];
\coordinate (A) at (intersection-1.center);
\coordinate (B) at (intersection-3.center);
  \path [name intersections={of=base and two}];
\coordinate (C) at (intersection-2.center);
\coordinate (D) at (intersection-4.center);

\coordinate (E) at (3cm,0);
\coordinate (F) at (0,3cm);
\coordinate (G) at (-3cm,0);
\coordinate (H) at (0,-3cm);
      

        \begin{scope}[even odd rule]
            \clip (A)--(0,0)--(F)--(0,5)--(5,5)--(A);
        \path[draw,black,dashed,line width=\lw] (0,0) ellipse \pone;
      \end{scope}
        \begin{scope}[even odd rule]
            \clip (C)--(0,2)--(H)--(5,-5)--(-5,-5)--(-5,0)--(C);
        \path[draw,black,dashed,line width=\lw] (0,0) ellipse \ptwo;
      \end{scope}
      \begin{scope}[even odd rule]
            \clip (G)--(E)--(5,5)--(-5,5)--(G);
        \path[draw,dashed,black,line width=\lw] (0,0) ellipse \pbase;
      \end{scope}

        \draw (A) circle (0.5ex)[fill=black]node[anchor=south west, black]{$12,13$};
        \draw (C) circle (0.5ex)[fill=black]node[anchor=-140, black]{$23,13$};

\shade [ball color=white,opacity=0.4] (0,0) circle (3cm);
      
        \begin{scope}[even odd rule]
            \clip (F)--(0,0)--(B)--(-5,5)--(5,5)--(F);
        \path[draw,black,line width=\lw] (0,0) ellipse \pone;
        \end{scope}
        \begin{scope}[even odd rule]
            \clip (H)--(0,0)--(D)--(5,-5)--(-5,-5)--(H);
        \path[draw,black,line width=\lw] (0,0) ellipse \ptwo;
        \end{scope} 
        \begin{scope}[even odd rule]
            \clip (G)--(E)--(5,-5)--(-5,-5)--(G);
        \path[draw,black,line width=\lw] (0,0) ellipse \pbase;
      \end{scope}

        \draw (B) circle (0.5ex)[fill=black]node[anchor=north east, black]{$13,12$};
        \draw (D) circle (0.5ex)[fill=black]node[anchor=north west, black]{$13,23$};
        \draw (E) circle (0.5ex)[fill=black]node[anchor=west, black]{$12,23$};
        \draw (G) circle (0.5ex)[fill=black]node[anchor=east, black]{$23,12$};
  \node[anchor=-135] at (45:3.2cm) {$v_1=0$};
  \node[anchor=135] at (-45:3.2cm) {$w_3=0$};
  \node[anchor=45] at (-135:3.2cm) {$v_3=0$};
  \node[anchor=-45] at (135:3.2cm) {$w_1=0$};
\end{tikzpicture}

\caption{\label{fig:pokeball} The totally nonnegative part of the complete flag variety in $\RR^3$. A vertex labeled $ab,cd$ denotes the cell $v_a\hspace*{-1pt}=\hspace*{-1pt}v_b\hspace*{-1pt}=\hspace*{-1pt}w_c\hspace*{-1pt}=\hspace*{-1pt}w_d\hspace*{-1pt}=\hspace*{-1pt}0$ (see \cref{ex:flag}).}
\end{figure}
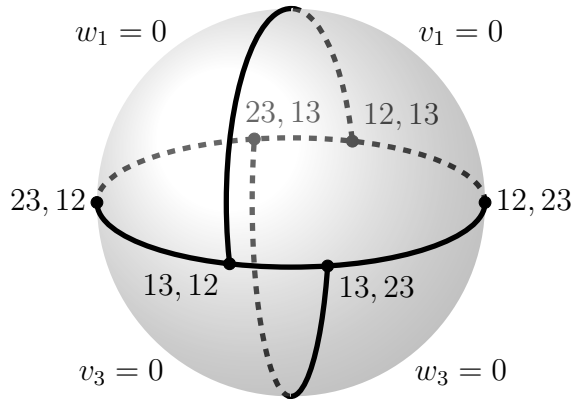

\subsection{Canonical basis}

We now assume that $G$ is simply laced. Let $Y$ (respectively, $X$) be the free abelian group of homomorphisms of algebraic groups $\RR^\ast\to T$ (respectively, $T\to \RR^\ast$), with standard pairing $\<\cdot,\cdot\>: Y\times X\to \Z$. Let $\DomWeights$ be the set of $\l\in X$ such that $\<\alpha_i^\vee,\l\>\geq 0$ for all $i\in\barI$, and for such $\l$ let $\supp(\l):=\{i\in\barI : \<\alpha_i^\vee,\l\>>0\}$ be its \emph{support}. 

For each $\l\in\DomWeights$, there is a unique irreducible $G$-module $\irrep_\l$ with highest weight $\l$. Lusztig~\Lus{3.1} defined a {\itshape canonical basis} $\CBl$ of $\irrep_\l$; see~\cite{LusCB} for more details. For $\line$ in the projective space $\PSl$, we write $\line>0$ (respectively, $\line\geq 0$) if the line $\line$ is spanned by a vector $x\in\irrep_\l$ whose coordinates in the canonical basis $\CBl$ are all positive (respectively, nonnegative). We denote
\[
\PSltp:=\{\line\in \PSl : \line>0\},\quad \PSltnn:=\{\line\in \PSl : \line\geq0\}.
\]

\def\openball{Q}
\def\closedball{\overline{\openball}}
\def\affinespan{R}
\def\flow{f}
\newcommand\act[2]{\flow(#1,#2)}

\subsection{Contractive flows}
We recall a topological lemma from~\cite{GKL} that we used to show that various spaces appearing in total positivity (including the totally nonnegative part of a type $A$ Grassmannian) are homeomorphic to closed balls.
\begin{defn}[{\cite[Definition~2.1]{GKL}}]\label{dfn:contract}
A map $\flow:\RR\times \RR^N\to\RR^N$ is called a \emph{contractive flow} if the following conditions are satisfied:
\begin{enumerate}[(1)]
\item\label{item:continuous} the map $\flow$ is continuous;
\item\label{item:action} for all $p\in\RR^N$ and $t_1,t_2\in\RR$, we have $\act 0 p=p$ and $\act{t_1+t_2}{p}=\act{t_1}{\act{t_2}{p}}$; and
\item\label{eq:contract} for all $p\neq 0$ and $t > 0$, we have $\|\act{t}{p}\| < \|p\|$.
\end{enumerate}
\end{defn}
Here $\|p\|$ denotes the Euclidean norm\footnote{In~\cite{GKL}, \cref{lem:top} was proved more generally for an arbitrary norm $\|\cdot\|$ on $\RR^N$.} of $p\in\RR^N$. For a subset $K\subset\RR^N$ and $t\in\RR$, let $\act{t}{K}$ denote $\{\act t p : p\in K\}$.

\begin{lem}[{\cite[Lemma~2.3]{GKL}}]\label{lem:top}
Let $\openball \subset \RR^N$ be a smooth embedded submanifold of dimension $d \leq N$, and $\flow:\RR\times \RR^N\to\RR^N$ a contractive flow. Suppose that $\openball$ is bounded and satisfies the condition
\begin{align}\label{eq:invariant}
\act{t}{\closedball} \subset \openball \quad \text{ for $t > 0$}.
\end{align}
Then the closure $\closedball$ is homeomorphic to a closed ball of dimension $d$, and $\closedball \setminus \openball$ is homeomorphic to a sphere of dimension $d-1$.
\end{lem}

\section{Proof of \texorpdfstring{\cref{main}}{Theorem 1}}\label{sec:proof-main-result}

\def\eig{\mu}
\def\M{\RR^N}
\def\p{p}
\def\H{\closedball}
\def\Hint{\openball}
\noindent We first prove \cref{main} in the case that $G$ is simply laced, in \cref{sec:simply-laced-case}. We will then deduce the result for all $G$ using the {\itshape folding} technique, in \cref{sec:general-case}. Our arguments employ the vector field $\tau := \sum_{i\in I}(e_i+f_i)\in\gfr$, and
\[
\expts:=\exp \left(t\tau\right)\text{ defined for }t\in\RR.
\]
By~\cite[Proposition~5.9(c)]{Lus}, we have $\expts\in\Gtp$ for $t>0$.

\begin{rmk}
In our earlier proof~\cite[\S 3]{GKL} that the totally nonnegative part of a type $A$ Grassmannian is homeomorphic to a closed ball, we worked with a closely related vector field, which is a cyclically symmetric (or loop group) analogue of $\tau$. We did so because exploiting cyclic symmetry appears to be necessary to show that the compactification of the space of electrical networks is homeomorphic to a ball, as we proved in~\cite[\S 6]{GKL}.
\end{rmk}

\subsection{Simply-laced case}\label{sec:simply-laced-case}

\def\emb{\psi}
Choose $\l\in \DomWeights$ such that $\supp(\l)=I\setminus J$. Then by~\LusGP{1.6}, for each $P\in\P^J$, there is a unique line $\line_P^\l$ in $\irrep_\l$ that is stable under the action of $P$ on $\irrep_\l$.  Moreover, the map $P\mapsto \line_P^\l$ is a smooth embedding of $\P^J$ into $\ProjSpace(\irrep_\l)$, which we denote by $\emb_J:\P^J\to\ProjSpace(\irrep_\l)$. Lusztig~\cite[Theorem~3.4]{LusGP} showed that we may select $\lambda$ so that the following hold:
\begin{enumerate}[(i)]
\item we have $P\in\Ptp^J$ if and only if $\line_P^\l>0$; and
\item we have $P\in\Ptnn^J$ if and only if $\line_P^\l\geq0$.
\end{enumerate}

Since $\expts\in\Gtp$ for $t>0$, by \cite[Theorem~5.6]{Lus} we have that $\tau$ acts on $\irrep_\l$ as a regular and semisimple endomorphism over $\mathbb{R}$. Let $v_0,v_1,\dots,v_N\in\irrep_\l$ be an eigenbasis of $\tau$ with eigenvalues $\eig_0 \ge \eig_1 \ge \cdots \ge \eig_N$ in $\RR$. We have a smooth open chart
\[
\phi:\M\hookrightarrow\PSl,
\]
which sends $\p=(\p_1,\dots,\p_N)\in\M$ to the line spanned by $v_0+\p_1v_1+\dots+\p_Nv_N$. By~\cite[Lemma~5.2(a)]{Lus} again, we have $\eig_k < \eig_0$ for $1\leq k\leq N$, and the image of $\phi$ contains $\PSltnn$. Consider the sets $\Hint\subset\H\subset\M$ defined by
\[
\Hint:=\phi^{-1}(\emb_J(\Ptp^J))\simeq\Ptp^J,\quad \H:=\phi^{-1}(\emb_J(\Ptnn^J))\simeq\Ptnn^J,
\]
so that $\H$ is compact and is the closure of $\Hint$.

Consider a map $\flow:\RR\times \RR^N\to\RR^N$ defined for $t\in\RR$ and $p=(p_1,\dots,p_N)\in\RR^N$ by
\begin{align}\label{coordinates}
\act t p:= \left(e^{t(\eig_1-\eig_0)}\p_1,\dots,e^{t(\eig_N-\eig_0)}\p_N\right).
\end{align}
We have $\phi(\act t p) = \expts\cdot\phi(p)$. We claim that $\flow$ is a contractive flow. Indeed, parts~\eqref{item:continuous} and~\eqref{item:action} of \cref{dfn:contract} hold for $\flow$. By~\eqref{coordinates}, we see that $\flow$ satisfies the property
\[
\| \act t p \| \le C^{-t} \| p \|\quad\text{ for all $t > 0$ and $p\in\RR^N$},
\]
where $C>1$ is the minimum of $e^{\eig_0-\eig_j}$ for $1 \le j \le N$. Therefore $f$ satisfies part~\eqref{eq:contract} of \cref{dfn:contract}, as claimed.

Lusztig~\cite[Proposition 2.2]{LusGP} showed that $\P^\barJ_{>0}$ is an open submanifold of $\P^\barJ$. Any $g\in\Gtp$ sends $\PSltnn$ inside $\PSltp$~\cite[Proposition~3.2]{Lus}, which implies that~\eqref{eq:invariant} holds. Therefore \cref{main} in the simply-laced case follows by applying \cref{lem:top}.

\begin{example}\label{ex:fixed_point}
Adopting the setup of \cref{ex:flag}, we have $\tau = \scalebox{0.72}{$\begin{bmatrix}0 & 1 & 0 \\ 1 & 0 & 1 \\ 0 & 1 & 0\end{bmatrix}$}$. The flow $f$ contracts $\P^\barJ_{\ge 0}$ onto the flag generated by the eigenvectors of $\tau$, ordered by decreasing eigenvalue. Its coordinates (see \cref{ex:flag}) are $v_1 = v_3 = w_1 = w_3 = \frac{1}{2+\sqrt{2}}, v_2 = w_2 = \frac{\sqrt{2}}{2+\sqrt{2}}$.
\end{example}

\subsection{General case, via folding}\label{sec:general-case}
For general $G$, by~\Lus{1.6} there exists an algebraic group $\dotG$ of simply-laced type with pinning $(\dotT,\dotBB^+,\dotBB^-,\dotx_i,\doty_i;i\in\dotI)$ and an automorphism $\sigma:\dotI\to\dotI$ that extends to an automorphism of $\dotG$ (also denoted $\sigma$), such that $G\cong \dotG^\sigma$. (Here $S^\sigma$ denotes the set of fixed points of $\sigma$ in $S$.) We have
\[
\Gtp=(\dotGtp)^\sigma,\quad \Gtnn=(\dotGtnn)^\sigma,\quad \Btp=(\dotBtp)^\sigma,\quad \Btnn=(\dotBtnn)^\sigma;
\]
see~\Lus{8.8}. It follows from~\eqref{eq:pi_J} that
\[
\Ptp^\barJ=(\dotPtp^\dotJ)^\sigma,\quad\Ptnn^\barJ=(\dotPtnn^\dotJ)^\sigma,
\]
where we identify $\barI$ with the set of equivalence classes of $\dotI$ modulo the action of $\sigma$, and $\dotJ$ is defined to be the preimage of $\barJ$ under the projection map $\dotI\to\barI$.

Note that $\sum_{i\in I}(e_i+f_i)\in\gfr$, so $\expts\in G$ for all $t\in\RR$. Thus $\expts$ acts on $\Ptnn^\barJ=(\dotPtnn^\dotJ)^\sigma$. The smooth embedding $\dotPtp^\dotJ \hookrightarrow \PSltp$ restricts to a smooth embedding $\Ptp^\barJ\hookrightarrow \PSltp$. Therefore we may apply \cref{lem:top} again to deduce that $\Ptnn^\barJ$ is homeomorphic to a closed ball, completing the proof of \cref{main} in the general case.

\bibliographystyle{alpha}
\bibliography{GP_bib}

\end{document}